\documentclass[12pt]{amsart}
\usepackage{amssymb}
\usepackage[all]{xy}
\usepackage{a4wide}
\newtheorem{thm}{Theorem}
\newtheorem{prop}[thm]{Proposition}

\newtheorem{conj}[thm]{Conjecture}
\newtheorem{rem}{Remark}

\newtheorem{lem}[thm]{Lemma}

\newenvironment{prf}{\begin{proof}[\bf{Proof}]}{\end{proof}}

\newcommand{\field}[1]{\mathbb{#1}}
\newcommand{\Q}{\field{Q}}
\newcommand{\C}{\field{C}}
\newcommand{\R}{\field{R}}

\newcommand{\N}{\field{N}}

\newcommand{\Z}{\field{Z}}

\newcommand{\ThetaS}{\Theta^{(S)}(\mathfrak{p},\,\mathfrak{a},\,\mathfrak{b})}
\newcommand{\OF}{\mathcal{O}_F}
\newcommand{\Ra}{\mathrm{R}_{\mathfrak{a}}}

\title[Computing Hilbert modular forms]{\bf Computing Hilbert modular forms over fields with nontrivial class group}
\author{Lassina Demb\'el\'e, Steve Donnelly}
\address{Institut f\"ur Experimentelle Mathematik\\ 
Ellernstrasse 29, 45326 Essen\\ Germany\\
e-mail: {\sf lassina.dembele@uni-due.de}}
\address{School of Mathematics and Statistics F07\\ University of Sydney NSW 2006\\ Sydney, Australia\\ email: {\sf donnelly@maths.usyd.edu.au}}

\begin{document}

\maketitle

\begin{abstract} In previous work, the first author developed an algorithm for the computation of Hilbert modular forms. In this paper, we extend this to all totally real number fields of even degree and nontrivial class group. 
Using the algorithm over $\Q(\sqrt{10})$ and $\Q(\sqrt{85})$ and their Hilbert class fields, we present some 
new instances of the conjectural Eichler-Shimura construction for totally real fields, and in particular find 
new examples of modular abelian varieties with everywhere good reduction.
\end{abstract}

\section*{\bf Introduction}
Let $F$ be a totally real number field of {\it even} degree. Let $\mathrm{B}$ be the unique quaternion algebra over $F$ which is ramified exactly at all infinite places. By the Jacquet-Langlands correspondence (Jacquet and Langlands~\cite[Chap. XVI]{jacqlang} and Gelbart~\cite{gelbart1}), computing Hilbert modular forms over $F$ amounts to computing automorphic forms on $\mathrm{B}$. In \cite{dembele1} and \cite{dembele2}, the first author presented an algorithm that exploits this correspondence. The algorithm adopts an alternative approach to the theory of Brandt matrices on $\mathrm{B}$ that is computationally more efficient than the classical one. Both papers considered only fields with narrow class number one.  One technical difficulty arising from nontrivial class groups is that ideals in $\mathrm{B}$ are no longer free $\mathcal{O}_F$-modules.  
This is now handled smoothly in the package for quaternion algebras over number fields included in the Magma computational  algebra system~\cite{magma} (version 2.14).  Our computations rely heavily on this package, in which algorithms from  \cite{voight1} and \cite{kirschmer} are implemented. There are not many explicit examples of Hilbert modular forms in the nontrivial class group case in the literature. Okada~ \cite{okada} provides several examples of such forms of level 1 and parallel weight 2 on the quadratic fields $\Q(\sqrt{257})$ and $\Q(\sqrt{401})$, computed using explicit trace formulae. However, it would be difficult to extend this approach to arbitrary totally real number fields, levels and weights.  In this paper we present a general algorithm that is practical for a large range of fields and levels.  This opens the possibility of experimenting systematically, especially over fields with nontrivial class group, and we hope this will shed new light on the theory of these objects.

The paper is organized as follows. Section~\ref{sec1} contains the necessary theoretical background. In section~\ref{sec2} we state the general algorithm, and describe some improvements to its implementation. Section~\ref{sec3} provides some numerical data over the real quadratic fields $\Q(\sqrt{10})$ and $\Q(\sqrt{85})$ and their Hilbert class fields. In section~\ref{sec4} we use this data to give new examples of the Eichler-Shimura construction over totally real number fields.

\medskip
\noindent
{\bf Acknowledgements.}  This project was started when the first author was still a PIMS postdoctoral fellow at the University of Calgary, and parts of it were written during his visit to the University of Sydney in August 2007. He would like to think both PIMS and the University of Calgary for their financial support, and the Department of Mathematics and Statistics of the University of Sydney for their hospitality. In particular, he would like to thank Anne and John Cannon for their invitation to visit the Magma group. He would also like to thank Clifton Cunningham for his constant support and encouragement in the early stage of the project. Finally, the authors would like to thank Fred Diamond and Noam Elkies for helpful email exchanges.

\section{\bf Theoretical background}\label{sec1}

Our aim is to compute spaces of Hilbert modular forms (as Hecke modules). By the Jacquet-Langlands correspondence, 
this is equivalent to computing spaces of automorphic forms $M_{\underline{k}}^{\mathrm{B}}(N)$ on some suitable 
quaternion algebra $\mathrm{B}$ (as Hecke modules).  In this section we define $M_{\underline{k}}^{\mathrm{B}}(N)$,
and then explain how to compute it in terms of spaces $M_{\underline{k}}({\mathrm{R}},N)$ of automorphic forms 
on quaternion orders $\mathrm{R} \subset \mathrm{B}$.  
A good reference for the material on Hilbert modular forms is Taylor~\cite{taylor1}. For the theory of Brandt matrices, 
we refer to Eichler~\cite{eichler1}, and also to \cite{dembele2} (which translates it into the adelic framework used here).  

We fix a totally real number field $F$ of {\it even} degree. We let $I$ be the set of all real embeddings of $F$; and for any $\tau\in I$, we denote the corresponding real embedding by $a\mapsto a^{\tau}$, $a\in F$. Also, we let $\mathcal{O}_F$ be the ring of integers of $F$, and fix an integral ideal $N$ of $F$. We let $\mathrm{B}$ be the
unique (up to isomorphism) quaternion algebra ramified at the infinite places of $F$. We fix a maximal order $\mathrm{R}$ 
of $\mathrm{B}$. We fix a Galois extension $K$ of $F$ contained in $\C$, which splits $\mathrm{B}$. We also fix 
an isomorphism $\mathrm{B}\otimes_F K\cong \mathbf{M}_2(K)^I$, and let $j:\,\mathrm{B}^\times\hookrightarrow
\mathbf{GL}_2(\C)^I$ be the resulting embedding. For each prime $\mathfrak{p}$ of $F$, we fix a local isomorphism $\mathrm{B}_{\mathfrak{p}}\cong\mathbf{M}_2(F_{\mathfrak{p}})$ such that
$\mathrm{R}_{\mathfrak{p}}=\mathbf{M}_2(\mathcal{O}_{F,\,\mathfrak{p}})$. Combining these local isomorphisms, one obtains
an isomorphism $\hat{\mathrm{B}}\cong\mathbf{M}_2(\hat{F})$ under which
$\hat{\mathrm{R}}\cong\mathbf{M}_2(\hat{\mathcal{O}}_F)$, where $\hat{F}$ and $\hat{\mathcal{O}}_F$ are the finite adeles 
of $F$ and $\mathcal{O}_F$ respectively. We define the compact open subgroup $U_0(N)$ of $\hat{\mathrm{R}}^\times$ by
$$U_0(N):=\left\{\begin{pmatrix}a&b\\ c&d\end{pmatrix}\in\mathbf{GL}_2(\hat{\mathcal{O}}_F):\,\,c\equiv 0(N)\right\}.$$

\medskip
Let $Cl_{\mathrm{B}}$ denote a complete set of representatives of all the right ideal classes of $\mathrm{R}$. 
Then $Cl_{\mathrm{B}}$ is in bijection with $\mathrm{B}^\times\backslash\hat{\mathrm{B}}^\times/\hat{\mathrm{R}}^\times.$
We choose a finite set of primes $S$ that generate the narrow class group $Cl_F^{+}$ and such that $\mathfrak{q}\nmid N$ for any $\mathfrak{q}\in S$. For any $\mathfrak{a}\in Cl_{\mathrm{B}}$, we let $\mathrm{R}_{\mathfrak{a}}$ be the left (maximal) order of $\mathfrak{a}$. Applying the strong approximation theorem, we may choose the representatives $\mathfrak{a}\in Cl_{\mathrm{B}}$ so that the prime divisors of $\mathrm{N}\mathfrak{a}$ lie in $S$. Then there are well-defined surjective reduction maps $\hat{\mathrm{R}}_{\mathfrak{a}}^\times\to\mathbf{GL}_2(\mathcal{O}_F/N)$, that all differ by conjugation in $\mathbf{GL}_2(\mathcal{O}_F/N)$. This gives a transitive action of each $\hat{\mathrm{R}}_{\mathfrak{a}}^\times$ on $\mathbf{P}^1(\mathcal{O}_F/N)$.

\medskip
Fix a vector $\underline{k}\in\Z^I$ such that $k_\tau\ge 2$ for all $\tau$, with all the components having the same parity. Set
$\underline{t}=(1,\,\ldots,\,1)$ and $\underline{m}=\underline{k}-2\underline{t}$, then choose $\underline{v}\in\Z^I$ such that each $v_\tau\ge 0$,
$v_\tau=0$ for some $\tau$, and $\underline{m}+2\underline{v}=\mu\underline{t}$ for some non-negative $\mu\in\Z$.  For every non-negative integer $a,\,b\in \Z$, we let $\mathbf{S}_{a,\,b}(\C)$ denote the right $\mathbf{M}_2(\C)$-module $\mathbf{Sym}^a(\C^2)$ (the $a^{th}$ symmetric power of the standard right $\mathbf{M}_2(\C)$-module $\C^2$) with the $\mathbf{M}_2(\C)$-action:
$$x\cdot m:=(\det m)^b x \mathrm{Sym}^a(m).$$ Then, we define the {\it weight} representation $L_{\underline{k}}$ by
$$L_{\underline{k}}=\bigotimes_{\tau\in I}\mathbf{S}_{m_\tau,\,v_\tau}(\C).$$

\medskip
The space of automorphic forms of level $N$ and weight $\underline{k}$ on $\mathrm{B}$ is defined as
%%%%\begin{eqnarray*}
$$
M_{\underline{k}}^{\mathrm{B}}(N):=\left\{f:\,\hat{\mathrm{B}}^\times/U_0(N) \to L_{\underline{k}}:\,\, f|\!|_{\underline{k}}\gamma=f,\quad \gamma\in\mathrm{B}^\times \right\},
$$
%%%%\end{eqnarray*}
where $f|\!|_{\underline{k}}\gamma(x):=f(\gamma x)\gamma.$

\begin{rem}\rm
By the Jacquet-Langlands correspondence~\cite[Chap. XVI]{jacqlang}, there is an isomorphism of Hecke modules between $M_{\underline{k}}^{\mathrm{B}}(N)$ and $M_{\underline{k}}(N)$, the space of Hilbert modular forms of weight $\underline{k}$ and level $N$ over $F$.  Therefore our task is now to compute $M_{\underline{k}}^{\mathrm{B}}(N)$ as a Hecke module.
\end{rem}

\noindent 
We define the space of automorphic forms of level $N$ and weight $\underline{k}$ on the order $\mathrm{R}_{\mathfrak{a}}$ by
$$M_{\underline{k}}(\mathrm{R}_\mathfrak{a},\,N):=\left\{f:\,\mathbf{P}^1(\mathcal{O}_F/N)
\to L_{\underline{k}}:\,\, f|\!|_{\underline{k}}\gamma=f,\quad \gamma\in \mathrm{R}_{\mathfrak{a}}^\times\right\}.$$ 

\medskip
For each $\mathfrak{a},\,\mathfrak{b}\in Cl_{\mathrm{B}}$ and any  prime $\mathfrak{p}\in F$, put
$$\ThetaS:=\mathrm{R}_{\mathfrak{a}}^\times\backslash\left\{u\in \mathfrak{a}\mathfrak{b}^{-1}:
   \,\,\frac{(\mathrm{nr}(u))}{\mathrm{nr}(\mathfrak{a})\mathrm{nr}(\mathfrak{b})^{-1}}=\mathfrak{p}\right\},$$ 
where $\mathrm{R}_{\mathfrak{a}}^\times$ acts by multiplication to the left. We define the linear map
\begin{eqnarray*}
T_{\mathfrak{a},\,\mathfrak{b}}(\mathfrak{p}):\,\,
M_{\underline{k}}(\mathrm{R}_{\mathfrak{b}},\,N) &\to& M_{\underline{k}}(\mathrm{R}_{\mathfrak{a}},\, N)\\
f&\mapsto&\sum_{u\in\Theta^{(S)}(\mathfrak{p},\,\mathfrak{a},\, \mathfrak{b})}f|\!|_{\underline{k}}u.
\end{eqnarray*}

\medskip
\begin{prop}\label{prop1} There is an isomorphism of Hecke modules
\begin{eqnarray*}
M_{\underline{k}}^{\mathrm{B}}(N)\to \bigoplus_{\mathfrak{a}\in Cl_{\mathrm{B}}} M_{\underline{k}}(\mathrm{R}_{\mathfrak{a}},\,N),
\end{eqnarray*} 
where the action of the Hecke operator $T(\mathfrak{p})$ on the right hand side is given by the collection of linear maps $(T_{\mathfrak{a},\,\mathfrak{b}}(\mathfrak{p}))$, $\mathfrak{a,\,b}\in Cl_{\mathrm{B}}$.
\end{prop}

\begin{prf} The Brandt matrices for a totally real number field are defined in Eichler~\cite{eichler1}, using the global language.  In that language, it is not hard to see that the action of the Hecke operator $T(\mathfrak{p})$ on the space $M_{\underline{k}}^{\mathrm{B}}(1)$ is determined by the collection of sets $\ThetaS$, $\mathfrak{a,\,b}\in Cl_{\mathrm{B}}$, for each prime $\mathfrak{p}$.  Let $G=\mathrm{Res}_{F/\Q}(\mathrm{B}^\times)$ and apply \cite[Theorem 1]{dembele3} to obtain Proposition~\ref{prop1}.  

Alternatively, one can observe that the proof of  \cite[Theorem 2]{dembele2} does not use the restriction on the class number of $F$. And so, Proposition~\ref{prop1} is simply a translation of \cite[Theorem 2]{dembele2} in global terms.
\end{prf}

\section{\bf Algorithmic issues}\label{sec2}
In the case of real quadratic fields, our algorithm has been discussed in \cite[sec. 2]{dembele1} and \cite[sec. 6]{dembele2}.  In this section we state the algorithm in full generality, for any totally real number field $F$ of even degree and any weight and level. We then discuss some implementation issues.

\medskip
As above, $\mathrm{B}$ will be the quaternion algebra ramified at the infinite places of $F$.
Our goal is to compute the space $M_{\underline{k}}^{\mathrm{B}}(N)$ of automorphic forms on $\mathrm{B}$
of weight $\underline{k}$ and level $N$, where $N$ is an integral ideal in $F$. 
First choose a (reasonable) bound $b\in\N,$ and a set of prime ideals $S$ not dividing $N$ that generate $Cl_F^{+}$.  
The more expensive tasks are done in an initial precomputation, since these depend only on $F$
and not on the weight or level.

\medskip
\noindent
{\bf Precomputation.}
\begin{enumerate}
\item 
Compute a maximal order $\mathrm{R}$ of $\mathrm{B}$.
\item 
Compute a complete set $Cl_{\mathrm{B}}$ of representatives $\mathfrak{a}$ for the right ideal classes of $\mathrm{R}$, 
chosen so that all prime factors of $\mathrm{N}\mathfrak{a}$ lie in $S$.
\item
For each representative $\mathfrak{a} \in Cl_{\mathrm{B}}$, compute its left order $\mathrm{R}_{\mathfrak{a}}$, 
and compute the unit group $\Gamma_{\mathfrak{a}}$ of $\mathrm{R}_{\mathfrak{a}}$.
\item 
Compute the sets $\ThetaS$, for all primes $\mathfrak{p}$ with $\mathrm{N}\mathfrak{p}\le b$ 
and all $\mathfrak{a,\,b}\in Cl_{\mathrm{B}}$.  (See Section~\ref{computing_Theta} for details.)
\end{enumerate}

\medskip \noindent
{\bf Algorithm.}
\begin{enumerate}
\item Compute splitting isomorphisms $\mathrm{R}_{\mathfrak{p}}^\times\cong\mathbf{GL}_2(\mathcal{O}_{F,\,\mathfrak{p}})$, for each prime $\mathfrak{p}\mid N$.
\item For each $\mathfrak{a}\in Cl_{\mathrm{B}}$, compute $M_{\underline{k}}(\mathrm{R}_{\mathfrak{a}},\,N)$ as a module of coinvariants
$$M_{\underline{k}}(\mathrm{R}_{\mathfrak{a}},\,N)=K[\mathbf{P}^1(\mathcal{O}_F/N)]\otimes L_{\underline{k}}/\langle x-\gamma x,\,\gamma\in\Gamma_{\mathfrak{a}}\rangle.$$
\item Compute the direct sum $$M_{\underline{k}}^{\mathrm{B}}(N)=\bigoplus_{\mathfrak{a}\in Cl_{\mathrm{B}}}M_{\underline{k}}(\mathrm{R}_{\mathfrak{a}},\,N).$$
\item For each prime $\mathfrak{p}$ of $F$ with $\mathrm{N}\mathfrak{p}\le b$, compute the families of linear maps $(T_{\mathfrak{a},\,\mathfrak{b}}(\mathfrak{p}))$.  (These determine the Hecke operator $T(\mathfrak{p})$ as a block matrix.)
\item Find a common basis of eigenvectors of $M_{\underline{k}}^{\mathrm{B}}(N)$ for the $T(\mathfrak{p})$.
\item If Step (5) does not completely diagonalize $M_{\underline{k}}^{\mathrm{B}}(N)$, increase $b$ and extend the precomputation, obtaining $\ThetaS$ for $\mathrm{N}\mathfrak{p}\le b$.  Then return to Step (4).
\end{enumerate}

\begin{rem}\rm
In practice, it is extremely rare that one resorts to Step (6) since very few Hecke operators $T(\mathfrak{p})$ 
are required to diagonalize the space $M_{\underline{k}}^{\mathrm{B}}(N)$. 
In the cases we tested, which included levels with norm as large as 5000, we never needed more than 10 primes.
\end{rem}

The steps in the main algorithm involve only local computations and linear algebra, whereas several steps in the 
precomputation involve lattice enumeration.  For a given field $F$, the precomputed data $\ThetaS$ can be 
re-used for all levels that are coprime to the primes in $S$, and all weights.  If one wishes to compute 
for all levels up to some large bound, one may simply choose the primes in $S$ to be larger than the bound.

\subsection{\bf Computing $\ThetaS$}
\label{computing_Theta}
Recall $\ThetaS \subset \Ra^{\times} \backslash \mathfrak{a}\mathfrak{b}^{-1}$ (as defined in Section~\ref{sec1}).
\begin{lem}
The correspondence $u \leftrightarrow u^{-1}\mathfrak{a}$ gives a bijection between $\ThetaS$ 
and the set of fractional right $\mathrm{R}$-ideals $\mathfrak{c} > \mathfrak{b}$ 
such that $\mathrm{nr}(\mathfrak{b}) = \mathrm{nr}(\mathfrak{c})\mathfrak{p}$
and $\mathfrak{c} \cong \mathfrak{a}$ as right $\mathrm{R}$-ideals.  
\end{lem}
\begin{prf}
Consider all fractional right $\mathrm{R}$-ideals $\mathfrak{c}$ that are isomorphic to $\mathfrak{a}$;
these are precisely the ideals $u^{-1}\mathfrak{a}$ for $u \in \mathrm{B}^{\times}.$  
Note that $u^{-1}\mathfrak{a} = v^{-1}\mathfrak{a}$ if and only if $v \in \Ra^{\times}u.$   
It is clear that $u^{-1}\mathfrak{a}$~contains~$\mathfrak{b}$ if and only if $u \in \mathfrak{a}\mathfrak{b}^{-1},$
and that $\mathrm{nr}(\mathfrak{b}) = \mathrm{nr}(u^{-1}\mathfrak{a})\mathfrak{p}$ 
if and only if $(\mathrm{nr}(u)) = \mathrm{nr}(\mathfrak{a})\mathrm{nr}(\mathfrak{b})^{-1}\mathfrak{p}.$
The lemma follows.
\end{prf}

\medskip \noindent {\bf Algorithm.}  
This computes $\ThetaS$ for all $\mathfrak{a} \in Cl_{\mathrm{B}}$, where $\mathfrak{p}$ and $\mathfrak{b}$ are fixed.
\begin{enumerate}
\item 
Compute the fractional right $\mathrm{R}$-ideals $\mathfrak{c} > \mathfrak{b}$ with 
$\mathrm{nr}(\mathfrak{b}) = \mathrm{nr}(\mathfrak{c})\mathfrak{p}$.
(The number of these is $\mathrm{N}\mathfrak{p} + 1.$)
\item 
For each such $\mathfrak{c}$, compute the representative ${\mathfrak{a}} \in Cl_{\mathrm{B}}$ and 
some $u \in \mathrm{B}$ such that $\mathfrak{c} = u^{-1} {\mathfrak{a}}$.  Append $u$ to $\ThetaS$. 
\end{enumerate}

\begin{rem}\rm
One sees that for each fixed $\mathfrak{p}$ and $\mathfrak{b},$
$$\sum_{\mathfrak{a} \in Cl_{\mathrm{B}}} \# \ThetaS = \mathrm{N}\mathfrak{p}+1\,,$$
however this fact is not used in the algorithm.
\end{rem}

Step (1) is a local computation; the ideals are obtained by pulling back local ideals 
under a splitting homomorphism $\mathrm{R}_{\mathfrak{p}}\cong\mathbf{M}_2(F_{\mathfrak{p}}).$   
Step (2) is the standard problem of isomorphism testing for right ideals, and we discuss 
an improvement to the standard algorithm for this in the next section; the complexity 
of each isomorphism test will not depend on $\mathfrak{p}.$

\subsection{\bf Lattice-based algorithms for definite quaternion algebras}
In this section, we let $\mathrm{B}$ be any definite quaternion algebra over a totally real number field $F$,
and let $\mathrm{R}$ be a (maximal or Eichler) order of $\mathrm{B}$. Two basic algorithmic problems are:  
\begin{enumerate}
\item 
to find an isomorphism between given right $\mathrm{R}$-ideals $\mathfrak{a}$ and $\mathfrak{b},$ and 
\item
to compute the unit group of $\mathrm{R}$ (modulo the unit group of $\OF$).
\end{enumerate}
The standard approach to both problems (as in \cite{voight1}) reduces them to the following.

\medskip \noindent {\bf General Problem:} 
{\it Let $L \cong \Z^n$ be a lattice contained in $\mathrm{B}$ (not necessarily of full rank).  
Given a totally positive element $\alpha \in F$, compute all $x \in L$ with $\mathrm{nr}(x) = \alpha$. }

\medskip 
For isomorphism testing, $L$ is the fractional ideal $\mathfrak{a}\mathfrak{b}^{-1}$ and 
$\alpha$ is a generator of $\mathrm{nr}(\mathfrak{a}\mathfrak{b}^{-1}).$ 
%%%%$\mathrm{nr}(\mathfrak{a})\mathrm{nr}(\mathfrak{b})^{-1}.$ 
(It suffices to consider a finite set of possibilities for $\alpha$.)
For computing units, $L$ is $\mathrm{R}$ (or possibly an $\OF$-submodule of $\mathrm{R}$ known to contain a unit),
and $\alpha$ is a unit of $\OF$.

One may solve the general problem by considering the positive definite quadratic form on $L$ 
given by $\mathrm{\mathrm{Tr}}(\mathrm{nr}(x)).$  (Note that its values are positive since 
$\mathrm{nr}(x)$ is a totally positive element of $F$, for all $0 \ne x \in \mathrm{B}$).
We capture all $x \in L$ with $\mathrm{nr}(x) = \alpha$ by enumerating all $x$ for which the 
quadratic form takes value $\mathrm{Tr}(\alpha)$, using standard lattice algorithms.
The drawback is that $\mathrm{Tr}(\alpha)$ might not be particularly small in relation to the determinant
of the lattice (even when $\alpha$ is a unit), in which case the lattice enumeration can be very time-consuming.  

We now present a variation which avoids this bottleneck.
For any nonzero $c \in F$, we may instead consider the lattice $cL \subset \mathrm{B}$, again under 
the positive definite quadratic form given by $\mathrm{Tr}(\mathrm{nr}(x))$.  
We clearly capture all $x \in L$ with $\mathrm{nr}(x) = \alpha$ by enumerating all $y \in cL$ with 
$\mathrm{Tr}(\mathrm{nr}(y)) = \mathrm{Tr}(c^2\alpha)$ and taking $x = y/c$.  In the special case 
that $c \in \Q$, this merely rescales the enumeration problem.  However, we will see that $c \in F$ 
may be chosen so that, in applications (1) and (2) above, one is looking for relatively short vectors
in the lattice.

Note that $\det(cL) = |\mathrm{N}(c)|^{\dim(L)/\deg F} \det(L)$.  Heuristically, as $c$ varies, the complexity 
of the enumeration process will be roughly proportional to the number of lattice elements with length 
up to the desired length.  Asymptotically this number equals
$$\frac{ \mathrm{Tr}(c^2\alpha)^{\dim(L)/2} }{ \det(cL) } \ = \ 
  \frac{ \mathrm{Tr}(c^2\alpha)^{\dim(L)/2} }{ |\mathrm{N}(c)|^{\dim(L)/\deg F} \det(L) } \ = \  
  \frac{ \mathrm{Tr}(c^2\alpha)^{\dim(L)/2} }{ \mathrm{N}(c^2\alpha)^{\dim(L)/{2\deg F}} } \ 
  \frac{ \mathrm{N}(\alpha)^{\dim(L)/{2\deg F}} }{\det(L)}\,,$$ 

Given that $\alpha$ is totally positive, $\mathrm{Tr}(c^2\alpha)/\mathrm{N}(c^2\alpha)^{1/\deg F}$ 
cannot be less than $\deg F$, and is close to $\deg F$ when all the real embeddings of $c^2\alpha$ 
lie close together.  It is straightforward to find $c \in \OF$ with this property, as follows.
First fix a $\Z$-module basis $\mathfrak{bas}(\OF)$ of $\OF$.

\medskip \noindent {\bf Algorithm.}\
Choose a large constant $C$.
\begin{enumerate} 
\item Calculate $r_i := C/{\sqrt{\sigma_i(\alpha)}}$ (note that the real embeddings of $\alpha$ are positive).
\item Represent the vector $(r_i)$ in terms of the basis $\mathfrak{bas}(\OF)$, then round the coordinates 
to integers, thus obtaining an element $c \in \OF$.
\end{enumerate} 

\begin{lem} Given any totally positive element $\alpha \in F$, and any $\epsilon > 0$, we can find $c \in F$ 
such that $\mathrm{Tr}(c^2\alpha)/\mathrm{N}(c^2\alpha)^{1/\deg F} < \deg F + \epsilon$.
\end{lem}

\begin{prf}
In the notation of the algorithm, we fix $\alpha$ and let $C \to \infty$, regarding $r_i \in\R$ 
and $c \in \OF$ as functions of $C$.  Since we use a fixed basis of $\OF$, $\sigma_i(c) - r_i$ is bounded by 
a constant independent of $C$.
Therefore as $C \to \infty$, $\sigma_i(c^2\alpha) = r_i^2\sigma_i(\alpha) + O(C) = C^2 + O(C) \,.$
This implies that for any $i$ and $j$, the ratio $\sigma_i(c^2\alpha)/\sigma_j(c^2\alpha) \to 1$ as $C \to \infty$,
and the lemma follows.
\end{prf}

The complexity of the enumeration thus depends on the ratio $\mathrm{N}(\alpha)^{\dim(L)/{2\deg F}}/{\det(L)}$.
In both the applications above, this ratio is small: in computing units, $\alpha$ is a unit, and in isomorphism
testing, $\alpha$ generates the fractional ideal $\mathrm{nr}(L)$ where $L = \mathfrak{a}\mathfrak{b}^{-1}$.

\section{\bf Examples of Hilbert modular forms}\label{sec3}

In this section we give some examples computed using our algorithm, which we have implemented in Magma 
(and which will be available in a future version of Magma). 

\subsection{\bf The quadratic field $\Q(\sqrt{85})$}\label{subsec31}
Let $F$ be the real quadratic field $\Q(\sqrt{85})$. The class number of $F$ is the same as its narrow class number: $h_F=h_F^+=2$. The maximal order in $F$ is $\mathcal{O}_F=\Z[\omega_{85}]$, where $\omega_{85}=\frac{1+\sqrt{85}}{2}$. Let $\mathrm{B}/F$ be the Hamilton quaternion algebra, i.e., the $F$-algebra given by 
$$\mathrm{B}:=F\oplus F i\oplus F j\oplus F k,\,\,\mbox{\rm with}\,\, i^2=-1, \,j^2=-1\,\,\mbox{\rm and}\,\, k=ij.$$ Since the prime $2$ is inert in $F$, the algebra $\mathrm{B}$ is ramified only at the two infinite places. Using Magma, we find that the class number of $\mathrm{B}$ is $8$. The Hecke module of Hilbert modular forms of level 1 and weight $(2,\,2)$ over $F$ is therefore an $8$-dimensional $\Q$-space, and it can be diagonalized by using the Hecke operator $T_2$.  There are two Eisenstein series and two Galois conjugacy classes of newforms. The eigenvalues of the Hecke operators for the first few primes are given in Table~\ref{table1} (only one eigenform in each Galois conjugacy class of newforms is listed).

\begin{table}
\begin{eqnarray*}
\begin{array}{|c|c||rrrr|}\hline
\mathrm{N}(\mathfrak{p})&\mathfrak{p}&f_1&f_2&f_3&f_4\\\hline\hline
3&(3, 2\omega_{85}) &4&-4& 2\sqrt{-1}& \beta \\ 
3&(3, 4+2\omega_{85}) &4 &-4& -2\sqrt{-1}& \beta\\
4&(2, 0) &5&5& 1& -\beta^3+3\\
5&(5, -1+2\omega_{85}) &6&-6& 0& -\beta^3+4\beta\\
7&(7, 2\omega_{85}) &8&-8& -2\sqrt{-1}& \beta^3-5\beta\\
7&(7, 12+2\omega_{85}) &8&-8& 2\sqrt{-1}& \beta^3-5\beta\\
17&(17, -1+2\omega_{85}) &18&-18& 0& 2\beta^3-14\beta\\
19&(19, 2+2\omega_{85})& 20&20& -4& 2\\
19&(19, 15+2\omega_{85}) &20&20& -4& 2\\\hline
%%%\mathrm{N}(\mathfrak{p})&\mathfrak{p}&f_1&f_2&f_3&f_4&f_5&f_6\\\hline\hline
%%%3&(3, 2\omega_{85}) &4&-4&-2\sqrt{-1}&2\sqrt{-1}&-\beta&\beta \\ 
%%%3&(3, 4+2\omega_{85}) &4 &-4&2\sqrt{-1}& -2\sqrt{-1}&-\beta&\beta\\
%%%4&(2, 0) &5&5& 1&1&-\beta^3+3&-\beta^3+3\\
%%%5&(5, -1+2\omega_{85}) &6&-6&0&0&\beta^3-4\beta&-\beta^3+4\beta\\
%%%7&(7, 2\omega_{85}) &8&-8&2\sqrt{-1}&-2\sqrt{-1}&-\beta^3+5\beta&\beta^3-5\beta\\
%%%7&(7, 12+2\omega_{85}) &8&-8&-2\sqrt{-1}&2\sqrt{-1}&-\beta^3+5\beta&\beta^3-5\beta\\
%%%17&(17, -1+2\omega_{85}) &18&-18&0&0&-2\beta^3+14\beta&2\beta^3-14\beta\\
%%%19&(19, 2+2\omega_{85})& 20&20&-4&-4&2&2\\
%%%19&(19, 15+2\omega_{85}) &20&20&-4&-4&2&2\\\hline
\end{array}
\end{eqnarray*}
\caption{Hilbert modular forms of level 1 and parallel weight 2 over $\Q(\sqrt{85})$. (The minimal polynomial of $\beta$ is $x^4-6x^2+2$).}
\label{table1}
\end{table}

\medskip
The Hilbert class field of $\Q(\sqrt{85})$ is $H := \Q(\sqrt{5},\, \sqrt{17}) = \Q(\alpha)$, where $\alpha^4 - 4\alpha^3 - 5\alpha^2 + 18\alpha - 1 = 0$.
The narrow class number of $H$ is $1$, and $\mathrm{B} \otimes_F H$ (the quaternion algebra over $H$ ramified at the four infinite places) has class number $4$. Thus the space of Hilbert modular forms of level $1$ and weight $(2,2)$ is 4-dimensional. The eigenvalues of the Hecke action for the first few primes are listed in Table~\ref{table2}. There is one Eisenstein series and two classes of newforms.  Elements of $\mathcal{O}_H$ are expressed in terms of the integral basis
$$ 1,\ \ \frac{1}{6}(\alpha^3 - 3\alpha^2 - 5\alpha + 10),\ \ 
     \frac{1}{6}(-\alpha^3 + 3\alpha^2 + 11\alpha - 10),\ \  
     \frac{1}{6}(-\alpha^3 + 14\alpha + 5), $$
which we use to write generators of the ideals in the table.
%%%\begin{eqnarray*}
%%%\beta_1&:=&1\\
%%%\beta_2&:=&\frac{1}{6}(\alpha^3 - 3\alpha^2 - 5\alpha + 10)\\
%%%\beta_3&:=&\frac{1}{6}(-\alpha^3 + 3\alpha^2 + 11\alpha - 10)\\
%%%\beta_4&:=&\frac{1}{6}(-\alpha^3 + 14\alpha + 5).
%%%\end{eqnarray*}

\begin{table}
\begin{eqnarray*}
\begin{array}{|c|c||rrr|}\hline
\mathrm{N}(\mathfrak{p})&\mathfrak{p}&f_1&f_2&f_3\\\hline\hline
4&[1, -1, 0, 1]&5&1&3+\beta'\\
4&[0, 2, -1, 1]&5&1&3+\beta'\\
9&[0, 1, -1, 0]&10&2&\beta'\\
9&[1, -1, -1, 0]&10&2&\beta'\\
19&[0, 1, 0, -1]&20&-4&2\\
19&[-1, 2, 0, 1]&20&-4&2\\
19&[1, -1, -1, 1]&20&-4&2\\
19&[-1, 2, -1, 1]&20&-4&2\\\hline
\end{array}
\end{eqnarray*}
\caption{Hilbert modular forms of level 1 and parallel weight 2 over the Hilbert class field $H$ of $\Q(\sqrt{85})$. (The minimal polynomial of $\beta'$ is $x^2+6x+2$).}
\label{table2}
\end{table}

%%%%%Finally, assume that we want to compute all the forms for norm up to 5000 on $\Q(\sqrt{85})$. Up to Galois conjugation, this will only require four precomputations with the primes $P_3$, $P_3'$, and $P_{37}$ (one of the prime above 37). For example, using the precomputation at $P_3$, we can compute all the forms for every prime level (up to Galois conjugation). As an illustration, we have listed some dimensions of spaces for all the prime levels of norm up to 20 in Table~\ref{table5}. (We see that they grow rather quickly).

We also computed some spaces over $\Q(\sqrt{85})$ with nontrivial level.  The dimensions of the spaces with prime level
of norm less than $100$ are given in Table~\ref{table5}.  (It suffices to consider just one prime in each pair
of conjugate primes, and for the precomputation we took $S = {(3,-1+\omega_{85})}$.)   For example, for level 
$\mathfrak{p} = (5,\sqrt{85})$, $M_2(\mathfrak{p})$ has dimension $20$, and the Hecke operator 
$T_{\mathfrak{q}}$ with $\mathfrak{q} = (7, 2\omega_{85})$ acting on $M_2(\mathfrak{p})$ has characteristic polynomial 
$$ (x - 8) (x + 8) (x^2 + 4)^2 (x^4 - 10x^2 + 18)^2 (x^6 + 28x^4 + 104x^2 + 100) .$$
Comparing this with the space $M_2((1))$ of level $1$, on which $T_{\mathfrak{q}}$ has characteristic polynomial
$$ (x - 8) (x + 8) (x^2 + 4) (x^4 - 10x^2 + 18) ,$$
one sees that the Hecke action on the subspace of newforms $M_2(\mathfrak{p})$ is irreducible, and the 
cuspidal oldform space embeds in $M_2(\mathfrak{p})$ under two degeneracy maps (as expected).

\begin{table}
\begin{eqnarray*}
\begin{array}{c|crrr|}\hline
\mathrm{N}(\mathfrak{p})&\dim M_2(\mathfrak{p})&\dim S_2(\mathfrak{p})&\dim S_2^{new}(\mathfrak{p})\\\hline\hline
3  &16&14&8\\
4  &24&22&16\\
5  &20&18&12\\
7  &32&30&24\\
17 &56&54&48\\
19 &68&66&60\\
23 &72&70&64\\
37 &124&122&116\\
59 &180&178&172\\
73 &232&230&224\\
89 &272&270&264\\
97 &304&302&296\\\hline\hline
%%%(3, -1+\omega_{85})&16&14&8\\
%%%(2)                &24&22&16\\
%%%(5, 2\omega_{85} -  1)&20&18&12\\
%%%(7, 2\omega_{85})     &32&30&24\\
%%%(17, 2\omega_{85} - 1)&56&54&48\\
%%%(19, 2\omega_{85} + 2)&68&66&60\\
%%%(23, 2\omega_{85} + 3)&72&70&64\\
%%%(37, 2\omega_{85} + 13)&124&122&116\\
%%%(59, 2\omega_{85} + 11)&180&178&172\\
%%%(73, 2\omega_{85} + 30)&232&230&224\\
%%%(89, 2\omega_{85} + 20)&272&270&264\\
%%%(97, 2\omega_{85} + 44)&304&302&296\\\hline\hline
\end{array}
\end{eqnarray*}
\caption{Dimensions of spaces of Hilbert modular forms over $\Q(\sqrt{85})$ with weight $(2,2)$ and prime level of norm less than $100$}
\label{table5}
\end{table}

\subsection{The quadratic field $\Q(\sqrt{10})$}\label{subsec32}
Let $F$ be the real quadratic field $\Q(\sqrt{10})$. The Hilbert class field of $F$ is 
$H:=\Q(\sqrt{2},\sqrt{5}) = \Q(\alpha)$, where the minimal polynomial of $\alpha$ is  $x^4 - 2x^3 - 5x^2 + 6x - 1$. 
The narrow class number of $H$ is $1$. We computed the space of Hilbert modular forms of level 1 and weight $(2,2)$ over $F$ and $H$, and the Hecke eigenvalues for the first few primes are listed in Table~\ref{table3} and Table~\ref{table4}  
(only one eigenform in each Galois conjugacy class of newforms is listed).
Elements of $\mathcal{O}_H$ are expressed in terms of the integral basis
$$ 1,\ \ \frac{1}{3}(2\alpha^3 - 3\alpha^2 - 10\alpha + 7),\ \  
   \frac{1}{3}(-2\alpha^3 + 3\alpha^2 + 13\alpha - 7),\ \  
   \frac{1}{3}(-\alpha^3 + 3\alpha^2 + 5\alpha - 8). $$
%%%\begin{eqnarray*}
%%%\beta_1&:=&1\\
%%%\beta_2&:=&\frac{1}{3}(2\alpha^3 - 3\alpha^2 - 10\alpha + 7)\\
%%%\beta_3&:=&\frac{1}{3}(-2\alpha^3 + 3\alpha^2 + 13\alpha - 7)\\
%%%\beta_4&:=&\frac{1}{3}(-\alpha^3 + 3\alpha^2 + 5\alpha - 8).
%%%\end{eqnarray*}
    
%Coefficients in integral basis: $$[  [0, 0, 1, 0], [1, 0, 1, -1], [0, 1, 0, 0], [-15, -44, -21, -26],  [-91, -123, -48, -97].$$
\begin{table}
\begin{eqnarray*}
\begin{array}{|c|c||r|r|r|}\hline
\mathrm{N}(\mathfrak{p})&\mathfrak{p}&f_1&f_2&f_3\\\hline\hline
2&(2, \omega_{40})&-3&3&-\sqrt{2}\\
3&( 3, \omega_{40} + 4)&-4&4&\sqrt{2}\\
3&( 3, \omega_{40} + 2)&-4&4&\sqrt{2}\\
5&( 5, \omega_{40})&-6&6&-2\sqrt{2}\\
13&( 13, \omega_{40} + 6)&-14&14&0\\
13&( 13, \omega_{40} + 7)&-14&14&0\\
31&( 31, \omega_{40} + 14)&32&32&4\\
31&( 31, \omega_{40} + 17)&32&32&4\\
37&( 37, \omega_{40} + 11)&-38&38&6\sqrt{2}\\
37&( 37, \omega_{40} + 26)&-38&38&6\sqrt{2}\\\hline
\end{array}
\end{eqnarray*}
\caption{Hilbert modular forms of level 1 and weight $(2,2)$ over $\Q(\sqrt{10})$.}
\label{table3}
\end{table}

\begin{table}
\begin{eqnarray*}
\begin{array}{|c|c||r|r|}\hline
\mathrm{N}(\mathfrak{p})&\mathfrak{p}&f_1&f_2\\\hline\hline
4&[0, 0, 1, 0]&5&-2\\
9&[1, 1, -1, 0]&10&-4\\
9&[0, 1, -1, 1]&10&-4\\
25&[1, -2, 0, 0]&26&-2\\
31&[1, 1, 1, -1]&32&4\\
31&[1, -1, -1, -1]&32&4\\
31&[1, 1, -1, 1]&32&4\\
31&[-3, 2, -1, 0]&32&4\\\hline
\end{array}
\end{eqnarray*}
\caption{Hilbert modular forms of level 1 and weight $(2,2)$ over the Hilbert class field $H$ of $\Q(\sqrt{10})$.}
\label{table4}
\end{table}

\section{\bf Examples of the Eichler-Shimura construction}\label{sec4}

In the study of Hilbert modular forms, the following conjecture is well-known. We refer to Shimura~\cite{shimura1} or Knapp~\cite{knapp1} for the classical case, and to Zhang~\cite{zhang1} and references therein for the number field case.

\begin{conj}[Eichler-Shimura]\label{conj1} Let $f$ be a Hilbert newform of level $N$ and parallel weight $2$ over a totally real field $F$. Let $K_f$ be the number field generated by the Fourier coefficients of $f$. Then there exists an abelian variety $A_f$ defined over $F$ such that $K_f\hookrightarrow\mathrm{End}(A_f)\otimes\Q$ and
$$L(A_f,\,s)=\prod_{\sigma\in\mathrm{Gal}(K_f/\Q)}L(f^\sigma,\,s),$$ where $f^\sigma$ is obtained by letting $\sigma$ act on the Fourier coefficients of $f$.
\end{conj}

In the classical setting, namely when $F=\Q$, this is a theorem known as the Eichler-Shimura construction. In general, many cases of the conjecture are also known. (See, for example, Zhang~\cite{zhang1} and references therein). In those cases the abelian variety $A_f$ is often constructed as a quotient of the Jacobian of some Shimura curve of level $N$.  The case where $[F:\,\Q]$ is even and $A_f$ has everywhere good reduction poses a different challenge. In this section, we provide new examples of such $A_f$.

\begin{rem} \rm
We refer back to the final paragraph of section~\ref{subsec31}.  The characteristic polynomials given there,
viewed in terms of Conjecture~\ref{conj1}, indicate that the newform part of $M_2(\mathfrak{p})$ corresponds to a 
simple abelian variety of dimension $6$. 
\end{rem}

\medskip
\subsection{The quadratic field $\Q(\sqrt{85})$} Keeping the notation of subsection~\ref{subsec31}, let $E/H$ 
be the elliptic curve with the following coefficients:

\begin{eqnarray*}
\begin{array}{|cccccc|}\hline
&a_1&a_2&a_3&a_4&a_6\\\hline\hline
E:&[1, 0, 0, 1]& [0, -1, 0, -1]& [0, 1, 1, 0]&[-5, -6, -1, 0]& [-8, -7, -3, 2]\\\hline
%%%%%E_2:& [1, 0, 1, 1]& [1, -1, -1, -1]& [1, 1, 1, 0]& [2, -4, -1, 0]& [-33, 18, -2, -1]\\\hline
\end{array}
\end{eqnarray*}

\medskip
\noindent
This curve has everywhere good reduction.  The restriction of scalars $A=\mathrm{Res}_{H/F}(E)$ 
is an abelian surface over $F$, also with everywhere good reduction.

\begin{rem}\rm
The $j$-invariant of $E$ is $64047678245 -12534349815\omega_{85} \in F$, and in fact $E$ is $H$-isomorphic to 
its conjugate under $\mathrm{Gal}(H/F)$.  Therefore $A$ is $H$-isomorphic to $E \oplus E$. Let $E^{\prime}$ denote one of the other two conjugates, which have $j$-invariant 
$51513328430+12534349815\omega_{85}$; there is an isogeny of degree $2$ from $E$ to $E^{\prime}$.
The restriction of scalars $\mathrm{Res}_{H/F}(E^{\prime})$ over $F$ also has good reduction everywhere;
it is $H$-isomorphic to $E^{\prime} \oplus E^{\prime}$, and is therefore isogenous to $A$.
\end{rem}

To establish the modularity of $E$ and $A$, we will apply the following result of Skinner and Wiles. 
Here we state the {\it nearly ordinary} assumption (Condition (iv)) in a slightly different way.

\begin{thm} \cite[Theorem A]{skinner1} \label{skinner-wiles} Let $F$ be a totally real abelian extension of $\Q$. Suppose that $p\ge 3$ is prime, and let $\rho:\,\mathrm{Gal}(\overline{F}/F)\longrightarrow\mathbf{GL}_2(\overline{\Q}_p)$ be a continuous, absolutely irreducible and totally odd representation unramified away from a finite set of places of $F$. Suppose that the reduction of $\rho$ is of the form $\bar{\rho}^{ss}=\chi_1\oplus\chi_2$, where $\chi_1$ and $\chi_2$ are characters, and suppose that:
\begin{itemize}
\item[{\it (i)}] the splitting field $F(\chi_1/\chi_2)$ of $\chi_1/\chi_2$ is abelian over $\Q$,
\item[{\it (ii)}]$(\chi_1/\chi_2)|_{D_v}\neq 1$ for each $v\mid p$,
\item[{\it (iii)}]$\rho|_{I_v}\cong\begin{pmatrix}\psi\epsilon_p^{k-1}&*\\ 0& 1\end{pmatrix}$ for each prime $v\mid p$,
\item[{\it (iv)}] $\det\rho=\psi\epsilon_p^{k-1}$, with $k\ge 2$ an integer, $\psi$ a character of finite order, and $\epsilon_p$ the $p$-adic cyclotomic character.
\end{itemize}
Then $\rho$ comes from a Hilbert modular from.
\end{thm}

\begin{prop}\label{prop2} $a)$  The elliptic curves $E$ is modular and corresponds to $f_2$ in Table~\ref{table2}.

\medskip
$b)$ The abelian surface $A$ is modular and corresponds to $f_3$ in Table~\ref{table1}.
\end{prop}

\begin{prf} $a)$ Let $\rho_{E,\,3}$ be the $3$-adic representation attached to $E$, and $\bar{\rho}_{E,\,3}$ the corresponding residual representation. Also, let $\mathfrak{p}$ be any prime above 3. Using Magma, we compute the torsion subgroup $E(H)_{tors}\cong\Z/2\oplus\Z/2$, and the trace of Frobenius $a_{\mathfrak{p}}(E)=2$. The latter implies that the representation $\rho_{E,\,3}$ is ordinary at $\mathfrak{p}$. By direct calculation, we find that $j(E)$ is the image of a
$H$-rational point on the modular curve $X_0(3)$:
$$j(E)=\frac{(\tau+27)(\tau+3)^3}{\tau},\,\mbox{\rm where}\,\,\tau=[2166, 527, -527, 1054].$$  This implies that $E$ has a Galois-stable subgroup of order 3, so the representation $\bar{\rho}_{E,\,3}$ is reducible. Since it is ordinary, there exist characters $\chi$, $\chi'$ unramified away from $\mathfrak{p}\mid 3$, with $\chi$ unramified at $\mathfrak{p}$, such that $\bar{\rho}_{E,\,3}^{ss}=\chi\oplus\chi'$ and $\chi\chi'=\epsilon_3$ is the mod $3$ cyclotomic character. The field $H(\chi/\chi')$ is clearly abelian. Therefore the representation $\rho_{E,\,3}$ satisfies the conditions of Skinner and Wiles, and $E$ is modular. Comparing traces of Frobenius with the eigenvalues given in Table~\ref{table2}, we see that the corresponding form is $f_2$. 

\medskip
$b)$ Let $\mathrm{BC}(f_3)$ be the base change from $F$ to $H$ of the newform $f_3$ in Table~\ref{table1}. Since the Hilbert class field extension $H/F$ is totally unramified, the form $\mathrm{BC}(f_3)$ has level 1 and trivial character. By comparing the Fourier coefficients at the split primes above 19, we see that $\mathrm{BC}(f_3)=f_2$ in Table~\ref{table2}. The result then follows from properties of restriction of scalars and base change.
\end{prf}

\begin{rem}\rm To find $E$, we reasoned as follows.  The eigenvalues of $f_2$ in Table~\ref{table2} suggest that the 
corresponding curve admits a $2$-isogeny.  The curve must have good reduction everywhere, and so must its conjugates;
if these are also modular, then they share the same $L$-series and are therefore isogenous to each other.  This would
mean the curve comes from an $H$-rational point on $X_0(2)$ whose $j$-invariant is integral.
Using a parametrisation of $X_0(2)$, we searched for such points.  We would like to thank Noam Elkies for suggesting this 
approach.
\end{rem}

\begin{rem}\rm If we assume Conjecture~\ref{conj1}, then there exists a modular abelian  surface $A$ over $H$ with real multiplication by $\Q(\sqrt{7})$ which corresponds to the form $f_3$ in Table~\ref{table2}. The restriction of scalars of $A$ from $H$ to $F$ is a modular abelian fourfold with real multiplication by $\Q(\beta)$ which corresponds to the form $f_4$ in Table~\ref{table1}.
\end{rem}

\medskip
\subsection{The quadratic field $\Q(\sqrt{10})$} Keeping the notation of subsection~\ref{subsec32}, let $E/H$ be the elliptic curve with the following coefficients:

\begin{eqnarray*}
\begin{array}{|cccccc|}\hline
&a_1&a_2&a_3&a_4&a_6\\\hline\hline
E:& [0, 0, 1, 0]& [1, 0, 1, -1]& [0, 1, 0, 0]& [-15, -44, -21, -26]&  [-91, -123, -48, -97]\\\hline
\end{array}
\end{eqnarray*}

\medskip
\noindent
This is an elliptic curve with everywhere good reduction over $H$. In contrast with the previous example, the four Galois conjugates have distinct $j$-invariants.  The restriction of scalars $A=\mathrm{Res}_{H/F}(E)$ is an abelian surface over $F$ with everywhere good reduction. 

\begin{prop}\label{prop3} The elliptic curve $E/H$ and the abelian surface $A/F$ are modular; $E$ corresponds to $f_2$ in Table~\ref{table4}, and $A$ corresponds to $f_3$ in Table~\ref{table3}.
\end{prop}

\begin{prf}
Let $\rho_{E,\,3}$ be the $3$-adic representation attached to $E$, and $\bar{\rho}_{E,\,3}$ its reduction modulo $3$. Then, $\bar{\rho}_{E,\,3}$ is reducible since $$j(E)=\frac{(\tau+27)(\tau+3)^3}{\tau},\,\mbox{\rm where}\,\,\tau=[5, 52, -18, -26].$$ As before, it is easy to see that $\rho_{E,\,3}$ satisfies the conditions of Skinner and Wiles.  So $E$ is modular, and hence $A$ is also modular. Comparing traces of Frobenius with Fourier coefficients, it is easy to see which forms in the tables
they correspond to.

Alternatively, we could consider the $7$-adic representation $\rho_{E,\,7}$. Its reduction mod $7$ is reducible since the point
$([16, 23, 9, 18]: [-157, -268, -119, -184]: [1, 0, 0, 0])$ is an $H$-rational point of order $7$ on $E$. Furthermore, for any prime $\mathfrak{p}\mid 7$, we have $a_{\mathfrak{p}}(E)=8$, and it is easy to see that $\rho_{E,\,7}$ satisfies the conditions of Skinner and Wiles.
\end{prf}

\begin{rem}\rm
It was shown by Kagawa~\cite[Theorem 3.2]{kagawa1} that there is no elliptic curve with everywhere good reduction over $\Q(\sqrt{10})$. Our results show that if we assume modularity in addition, there is only one such simple abelian variety: an abelian surface with real multiplication by $\Z[\sqrt{2}]$. 
\end{rem}

\begin{rem}\rm To find $E$, we were again assisted by the eigenvalues of the corresponding form $f_2$ in Table~\ref{table4}, which suggest that $E$ has an $H$-rational point of order 14. The modular curve $X_0(14)/\Q$ is an elliptic curve (14A1 in Cremona's table), which (using Magma) was found to have rank $1$ over $H$ and also rank $1$ over $\Q(\sqrt{10})$; this enabled us to obtain a point of infinite order simply by finding a $\Q$-rational point on the quadratic twist by $\sqrt{10}$. We considered curves corresponding to points of small height in $X_0(14)(H)$, and twists of these curves, until we found one with good reduction everywhere.
\end{rem}


\begin{thebibliography}{99}
\bibitem{magma} Wieb Bosma, John Cannon, and Catherine Playoust. The Magma algebra system. I. The user language. {\it J. Symbolic Comput.}, {\bf 24}(3-4): 235-265, 1997.

\bibitem{dembele1} L. Demb\'el\'e, Explicit computations of Hilbert modular forms on $\mathbb{Q}(\sqrt{5})$. {\it Experiment. Math.} {\bf 14} (2005), no. 4, 457--466.

\bibitem{dembele2} L. Demb\'el\'e, Quaternionic $M$-symbols, Brandt matrices and Hilbert modular forms.  {\it Math. Comp.} {\bf 76}, no 258, (2007), 1039--1057.

\bibitem{dembele3} L. Demb\'el\'e, On the computation of algebraic modular forms. (submitted).

\bibitem{eichler1} M. Eichler, On theta functions of real algebraic number fields. {\it Acta Arith.} {\bf 33} (1977), no 3, 269--292.

\bibitem{gelbart1} S. Gelbart, Automorphic forms on adele groups. {\it Annals of Maths. Studies} {\bf 83}, Princeton Univ. Press, 1975.

\bibitem{jacqlang} H. Jacquet and R. P. Langlands, Automorphic forms on GL(2). {\it Lectures Notes in Math.}, vol. {\bf 114}, Springer-Verlag, Berlin and New York, 1970.

\bibitem{kagawa1}T. Kagawa, Elliptic curves with everywhere good reduction over real quadratic fields. Ph. D Thesis, Waseda University, 1998.

\bibitem{kirschmer} M. Kirschmer, Konstruktive Idealtheorie in Quaternionenalgebren.  Diplom Thesis, Universit\"at Ulm, 2005.

\bibitem{knapp1} A. W. Knapp, {\it Elliptic curves}. Mathematical Notes, {\bf 40}. Princeton University Press, Princeton, NJ, 1992. xvi+427 pp.

\bibitem{okada} Okada, Kaoru; Hecke eigenvalues for real quadratic fields. {\it Experiment. Math.} {\bf 11} (2002), no. 3, 407--426.

\bibitem{shimura1} G. Shimura, {\it Introduction to the arithmetic theory of automorphic functions.} Kan\^o Memorial Lectures, No. 1. Publications of the Mathematical Society of Japan, No. 11. Iwanami Shoten, Publishers, Tokyo; Princeton University Press, Princeton, N.J., 1971. xiv+267 pp.

\bibitem{skinner1} Skinner, C. M.; Wiles, A. J. Residually reducible representations and modular forms. {\it Inst. Hautes \'Etudes Sci. Publ. Math.} No. {\bf 89} (1999), 5--126 (2000).

\bibitem{taylor1} R. Taylor, On Galois representations associated to Hilbert modular forms. {\it Invent. Math}. {\bf 98} (1989), no. 2, 265--280. 

\bibitem{voight1} J. Voight, Quadratic forms and quaternion algebras: Algorithms and arithmetic. Ph. D thesis, University of Berkeley, 2005.

\bibitem{zhang1}  Zhang, Shouwu; Heights of Heegner points on Shimura curves. {\it Ann. of Math.} (2) {\bf 153} (2001), no. 1, 27--147. 

\end{thebibliography}
\end{document}